\documentstyle[11pt, amsfonts]{amsart}
\begin{document}
\def\R{{\mathbb R}}
\def\Rad{{\mathcal R}}
\def\Z{{\mathbb Z}}
\def\C{{\mathbb C}}
\newcommand{\trace}{\rm trace}
\newcommand{\Ex}{{\mathbb{E}}}
\newcommand{\Prob}{{\mathbb{P}}}
\newcommand{\E}{{\cal E}}
\newcommand{\F}{{\cal F}}
\newtheorem{df}{Definition}
\newtheorem{theorem}{Theorem}
\newtheorem{lemma}{Lemma}
\newtheorem{pr}{Proposition}
\newtheorem{co}{Corollary}
\newtheorem{pb}{Problem}
\def\n{\nu}
\def\sign{\mbox{ sign }}
\def\a{\alpha}
\def\N{{\mathbb N}}
\def\A{{\cal A}}
\def\L{{\cal L}}
\def\X{{\cal X}}
\def\F{{\cal F}}
\def\c{\bar{c}}
\def\v{\nu}
\def\d{\delta}
\def\diam{\mbox{\rm dim}}
\def\vol{\mbox{\rm Vol}}
\def\b{\beta}
\def\t{\theta}
\def\l{\lambda}
\def\e{\varepsilon}
\def\colon{{:}\;}
\def\pf{\noindent {\bf Proof :  \  }}
\def\endpf{ \begin{flushright}
$ \Box $ \\
\end{flushright}}
\title[Inequalities for the Radon transform]
{Inequalities for the derivatives of the Radon transform on convex bodies}

\author{Wyatt Gregory}
\address[Wyatt Gregory]{Department of Mathematics, University of Missouri, Columbia, MO 65211, USA}
\email{wtgmtf@@mail.missouri.edu}

\author{Alexander Koldobsky}
\address[Alexander Koldobsky]{Department of Mathematics, University of Missouri, Columbia, MO 65211, USA}
\email{koldobskiya@@missouri.edu}

\thanks{The second named author was supported by the U.S. National Science Foundation Grant DMS-1700036.}
\date{}

\begin{abstract}  It was proved in \cite{K-2014} that the sup-norm of the Radon transform of an arbitrary probability density on an origin-symmetric convex body of volume 1 is bounded from below by a positive constant depending only on the dimension. In this note we extend this result to the derivatives of the Radon transform. We also prove a comparison theorem for these derivatives.
\end{abstract}  
\maketitle
\section{Introduction} Let $K$ be an origin-symmetric convex body of volume 1 in $\R^n,$ and let $f$ be any non-negative measurable function on $K$ with $\int_Kf=1.$ Does there exist a constant $c_n$ depending only on $n$ so that for any such $K$ and $f$ there exists a direction $\xi\in S^{n-1}$ with $\int_{K\cap \xi^\bot}f \ge c_n?$ Here $\xi^\bot=\{x\in \R^n: (x,\xi)=0\}$ is the central hyperplane perpendicular to $\xi,$ and integration is with respect to Lebesgue measure on $\xi^\bot.$ It was proved in \cite{K-2014} that, in spite of the generality of the question, the answer to this question is positive, and one can take $c_n>\frac 1{2\sqrt{n}}.$ In \cite{CGL} this result was extended to non-symmetric bodies $K.$ Moreover, it was shown in \cite{KK} that this estimate is optimal up to a logarithmic term, and the logarithmic term was removed in \cite{KLi}, so, finally, $c_n\sim \frac 1{\sqrt{n}}.$ We write $a\sim b$ if there exist absolute constants $c_1,c_2>0$ such that $c_1a\le b \le c_2a.$  

Note that the same question for volume, where $f\equiv 1,$ is the matter of the slicing problem of Bourgain \cite{Bo1, Bo2}. In this case the best known result $c_n>cn^{-\frac 14},$ where $c>0$ is an absolute constant, is due to Klartag \cite{Kla} who removed a logarithmic term from an earlier result of Bourgain \cite{Bo3}. The constant does not depend on the dimension for several classes of bodies $K.$ For example, it was proved in \cite{K2} that if $K$ belongs to the class of unconditional convex bodies, the constant $c_n=\frac 1{2e}$ works for all functions $f.$ The same happens for intersection bodies \cite{K-2012}, and for the unit balls of subspaces of $L_p,\ p>2$ where the constant is of the order $p^{-1/2}$ \cite{KP}. 

Denote by
$$Rf(\xi,t) = \int_{K\cap \{x\in \R^n:(x,\xi)=t\}} f(x) dx,\qquad \xi\in S^{n-1},\ t\in \R$$
the Radon transform of $f.$ The result described above means that the sup-norm of the Radon transform of a probability density on an origin-symmetric convex body in $\R^n$ is bounded from below by a positive constant depending only on the dimension $n.$

In this note we prove a similar estimate for the derivatives of the Radon transform. For an origin-symmetric infinitely smooth convex body $K$ in $\R^n,$ an infinitely differentiable function $f$ on $K$, and $q\in \C,\ \Re q>-1,$ we denote the fractional derivative of the order $q$ of the Radon transform of $f$ by

$$(Rf(\xi,t))_t^{(q)}(0)=\left(\int_{K\cap \{x:\ (x,\xi)=t\}} f(x) dx\right)_t^{(q)}(0).$$

The estimate that we prove is as follows.

\begin{theorem} \label{fract-der}
There exists an absolute constant $c>0$ so that for any infinitely smooth origin-symmetric convex body $K$ of volume 1 in $\R^n,$ any even infinitely smooth probability density $f$ on $K,$ and any $q\in \R,\ 0\le q\le n-2,$ which is not an odd integer, there exists a direction $\xi\in S^{n-1}$ so that
$$\left(c \frac {q+1}{\sqrt{n\log^3(\frac{ne}{q+1})}}\right)^{q+1}\le  \frac 1{\cos(\frac{\pi q}2)} (Rf(\xi,t))_t^{(q)}(0).$$
 \end{theorem}
 
 Compare this with (\ref{eucl-norm}) where it is shown that for $f\equiv 1$ and $K=B_n$ the Euclidean ball of volume 1,
 the constant in the left-hand side does not depend on $n,$ and is greater or equal to $(c(q+1))^{\frac{q+1}2},$ where $c$ is an absolute constant. We deduce Theorem \ref{fract-der} from a more general result.
 
 \begin{theorem} \label{slicing} Suppose $K$ is an infinitely smooth origin-symmetric convex body in $\R^n$, $f$ is a non-negative even infinitely smooth function on $K$ and $-1<q<n-1$ is not an odd integer. Then
$$\int_K f \le \frac {n}{(n-q-1)\ 2^q \pi^{\frac{q-1}2} \Gamma(\frac{q+1}2)} |K|^{\frac{q+1}n} \left(d_{\rm ovr}(K,L^n_{-1-q})\right)^{q+1}$$
$$\times \max_{\xi\in S^{n-1}} \frac 1{\cos(\frac{\pi q}2)} (Rf(\xi,t))_t^{(q)}(0).$$
\end{theorem}

Here $|K|$ stands for volume of proper dimension, and the outer volume ratio distance from  $K$ to the class $L_{-1-q}^n$ of the unit balls of $n$-dimensional spaces that embed in $L_{-1-q}$ (see definition below) is defined by
$$d_{{\rm {ovr}}}(K,L_{-1-q}^n) = \inf \left\{ \left( \frac {|D|}{|K|}\right)^{1/n}:\ K\subset D,\ D\in L_{-1-q}^n \right\}.$$

For some classes of bodies $K$ the constant in the left-hand side of the estimate of Theorem \ref{fract-der} does not depend on the dimension, and behaves like in the Euclidean case (formula (\ref{eucl-norm})). This follows from estimates for the distance $d_{\rm ovr}(K,L^n_{-1-q}).$ Indeed, this distance is equal to 1 when $K$ is an intersection body. It is less or equal to $e$ for unconditional convex bodies $K,$ and it is less than $c\sqrt{p}$ for the unit balls of subspaces of $L_p$ with $p>2,$ where $c$ is an absolute constant; see remarks at the end of the paper. Also note
that the unit balls of subspaces of $L_p$ with $0<p\le 2$ are intersection bodies, so the distance is 1; see \cite{K5}.
 
 Our next result is related to the Busemann-Petty problem \cite{BP} which asks the following question. Let $K,L$ be origin-symmetric convex bodies in $\R^n,$ and suppose that the $(n-1)$-dimensional volume of every central hyperplane section of $K$ is smaller than the same for $L,$ i.e. $|K\cap \xi^\bot|\le |L\cap\xi^\bot|$ for every $\xi\in S^{n-1}.$ Does it necessarily follow that the $n$-dimensional volume of $K$ is smaller than the volume of $L,$ i.e. $|K|\le |L|?$ The answer is affirmative if the dimension $n\le 4,$ and it is negative when $n\ge 5;$ see \cite{G3, K1} for the solution and its history.
 
 It was proved in \cite{K-1999} (see also \cite[Theorem 5.12]{K1}) that the answer to the Busemann-Petty problem becomes affirmative if one compares the derivatives of the parallel section function of high enough orders. Namely, denote by
$$A_{K,\xi}(t) = R(\chi_K)(\xi,t)= |K\cap \{x\in \R^n:\ (x,\xi)=t\}|,\quad t\in \R$$
the parallel section function of $K$ in the direction $\xi.$ If $K,L$ are infinitely smooth origin-symmetric convex bodies in $\R^n,\ n\ge 4$, $q\in [n-4,n-1)$ is not an odd integer, and for every $\xi\in S^{n-1}$ the fractional derivatives of the order $q$ of the parallel section functions at zero satisfy
$$\frac 1{\cos(\frac{\pi q}2)}A_{K,\xi}^{(q)}(0)\le \frac 1{\cos(\frac{\pi q}2)}A_{L,\xi}^{(q)}(0),$$
then $|K|\le |L|.$ For $-1<q<n-4$ this is no longer true.

Another generalization of the Busemann-Petty problem, known as the isomorphic
Busemann-Petty problem, asks whether the inequality for volumes holds up to an absolute constant.  Does there exist an absolute constant $C$ so that for any dimension $n$ and any origin-symmetric
convex bodies $K,L$ in $\R^n$ satisfying
$|K\cap \xi^\bot|\le |L\cap \xi^\bot|$ for all $\xi\in S^{n-1},$
we have $ |K|\le C |L| ?$ This question is equivalent to the slicing problem of Bourgain mentioned above.

Zvavitch \cite{Zv} considered an extension of the Busemann-Petty problem to general Radon transforms, as follows. Suppose that $K,L$ are origin-symmetric convex bodies in $\R^n$, and $f$ is an even continuous strictly positive function on $\R^n.$ Suppose that $R(f\vert_K)(\xi,t)(0)\le R(f\vert_L)(\xi,t)(0)$ for every $\xi\in S^{n-1},$ where $f\vert_K$ is the restriction of to $K.$
Does it necessarily follow that $|K|\le |L|?$ Isomorphic versions of this result were proved in \cite{KZ, KPZv}.
 
 In this note we generalize these results to general Radon transforms as follows.
 
\begin{theorem} \label{bp} Let $K,L$ be infinitely smooth origin-symmetric convex bodies in $\R^n,$ $f,g$ non-negative infinitely differentiable functions on $K$ and $L$, respectively, $\|g\|_{\infty}=g(0)=1,$  and $q\in (-1,n-1)$ is not an odd integer. If for every $\xi\in S^{n-1}$
 $$\frac 1{\cos(\frac{\pi q}2)}(Rf(\xi,t))_t^{(q)}(0)\le \frac 1{\cos(\frac{\pi q}2)}(Rg(\xi,t))_t^{(q)}(0),$$
 then
 $$\int_K f \le \frac n{n-q-1} \left(d_{\rm ovr}(K,L^n_{-1-q})\right)^{q+1}\left(\int_L g\right)^{\frac {n-q-1}{n}} |K|^{\frac {q+1}n} .$$
 \end{theorem}
 
 In the case $q=0$ this result was proved in \cite{KPZv}.

\section{Notation and auxiliary facts} 

A closed bounded set $K$ in $\R^n$ is called a {\it star body}  if 
every straight line passing through the origin crosses the boundary of $K$ 
at exactly two points, the origin is an interior point of $K,$
and the {\it Minkowski functional} 
of $K$ defined by 
$\|x\|_K = \min\{a\ge 0:\ x\in aK\}$
is a continuous function on $\R^n.$ If $x\in S^{n-1},$ then $\|x\|_K^{-1}=r_K(x)$
is the radius of $K$ in the direction of $x.$ A star body $K$ is {\it origin-symmetric} if $K=-K.$
A star body $K$ is called {\it convex} if for any $x,y\in K$ and every $0<\lambda<1,$
$\lambda x+ (1-\lambda)y\in K.$

The polar formula for the volume of a star body is as follows:
\begin{equation} \label{polar}
|K|=\frac1n \int_{S^{n-1}} \|x\|_K^{-n} dx.
\end{equation}
A star body $K$ in $\R^n$ is called {\it infinitely smooth} if $\|x\|_K\in C^{\infty}(S^{n-1}).$

We use the techniques of the Fourier approach
to sections of convex bodies that has recently been developed; see \cite{K1}. 
The Fourier transform of a
distribution $f$ is defined by $\langle\hat{f}, \phi\rangle= \langle f, \hat{\phi} \rangle$ for
every test function $\phi$ from the Schwartz space $ \mathcal{S}$ of rapidly decreasing infinitely
differentiable functions on $\R^n$. For any even distribution $f$, we have $(\hat{f})^\wedge
= (2\pi)^n f$.

If $K$ is a star body  and $0<p<n,$
then $\|\cdot\|_K^{-p}$  is a locally integrable function on $\R^n$ and represents a distribution. 
Suppose that $K$ is infinitely smooth, i.e. $\|\cdot\|_K\in C^\infty(S^{n-1})$ is an infinitely differentiable 
function on the sphere. Then by \cite[Lemma 3.16]{K1}, the Fourier transform of $\|\cdot\|_K^{-p}$  
is an extension of some function $g\in C^\infty(S^{n-1})$ to a homogeneous function of degree
$-n+p$ on $\R^n.$ When we write $\left(\|\cdot\|_K^{-p}\right)^\wedge(\xi),$ we mean $g(\xi),\ \xi \in S^{n-1}.$
If $K,L$ are infinitely smooth star bodies, the following spherical version of Parseval's
formula was proved in \cite{K-1999} (see \cite[Lemma 3.22]{K1}):  for any $p\in (-n,0)$
\begin{equation}\label{parseval}
\int_{S^{n-1}} \left(\|\cdot\|_K^{-p}\right)^\wedge(\xi) \left(\|\cdot\|_L^{-n+p}\right)^\wedge(\xi) =
(2\pi)^n \int_{S^{n-1}} \|x\|_K^{-p} \|x\|_L^{-n+p}\ dx.
\end{equation}

A distribution is called {\it positive definite} if its Fourier transform is a positive distribution in
the sense that $\langle \hat{f},\phi \rangle \ge 0$ for every non-negative test function $\phi.$

In this note we extensively use fractional derivatives. In general, the
fractional derivative of the order $q$ of a test function $\phi\in {\mathcal S}(\R)$ 
is defined as the convolution of $\phi$ with $t_+^{-1-q}/\Gamma(-q):$
$$\phi^{(q)}(x) = \langle \frac{t_+^{-1-q}}{\Gamma(-q)}, \phi(x-t) \rangle .$$
We need this definition only with $x=0.$ This allows us to replace the test function
$\phi$ by any function differentiable up to a certain order in a neighborhood of zero.
 
Let $m\in \N\cup \{0\}$ and suppose that $h$ is a continuous
function on $\R$ that is $m$ times continuously differentiable 
in some neighborhood of zero. 
For $q\in \C,\ -1<\Re(q)<m,\ q\neq 0,1,...,m-1,$
the {\it fractional derivative} of the order $q$ of the function $h$ at zero
is defined as the action of the distribution $t_+^{-1-q}/\Gamma(-q)$ on the 
function $h,$ as follows: 
$$h^{(q)}(0) = \frac{1}{\Gamma(-q)}
\int_0^1 t^{-1-q}\big(h(t)-h(0)-...-h^{(m-1)}(0)
\frac{t^{m-1}}{(m-1)!}\big)\ dt +$$
\begin{equation}\label{fracder}
\frac{1}{\Gamma(-q)}\int_1^\infty t^{-1-q}h(t) dt +
\frac{1}{\Gamma(-q)} \sum_{k=0}^{m-1} \frac{h^{(k)}(0)}{k!(k-q)}.
\end{equation}
It is easy to see that for a fixed $q$ the definition does not depend on
the choice of $m>\Re(q),$ as long as $f$ is $m$ times continuously
differentiable. Note that without dividing by $\Gamma(-q)$ the expression 
for the fractional derivative represents an analytic function in the 
domain $\{q\in \C:\ \Re(q)>-1\}$ not including integers, and has simple 
poles at integers. The function $\Gamma(-q)$ is analytic 
in the same domain and also has simple poles at non-negative integers,
so after the division we get an analytic function in the whole domain
$\{q\in \C:\ m>\Re(q)>-1\},$ which also defines fractional derivatives 
of integer orders.  Moreover, computing the limit as $q\to k,$
where $k$ is a non-negative integer, we see that
the fractional derivatives of integer orders coincide with usual 
derivatives up to a sign (when we compute the limit the first two summands
in the right-hand side of (\ref{fracder}) converge to zero, since $\Gamma(-q)\to \infty,$
and the limit in the third summand can be computed using the property $\Gamma(x+1)=x\Gamma(x)$
of the $\Gamma$-function):
$$h^{(k)}(0) = (-1)^k \frac{d^k}{d t^k} 
h(t)\vert_{t=0}.$$
If $h$ is an even function its derivatives of odd orders at the origin 
are equal to zero and, for $m-2<\Re(q)<m,$ the expression (\ref{fracder}) 
turns into
\begin{equation}\label{fracdereven}
h^{(q)}(0)
={1\over \Gamma(-q)}\int^\infty_0 t^{-q-1} 
\left(h(t) - \sum_{j=0}^{ (m-2)/2} {t^{2j}\over (2j)!} 
h^{(2j)}(0)\right) d t.
\end{equation}
We also note that if $-1<q<0$ then
\begin{equation} \label{fracder-1}
h^{(q)}(0) = \frac{1}{\Gamma(-q)}\int_0^\infty t^{-1-q} h(t)\ dt.
\end{equation}

We need a simple fact that can be found in \cite[Lemma 3.14]{K1}.
\begin{lemma} \label{lemma3.14} Let $0<q<1.$ For every even test function $\phi$ and every fixed vector $\theta\in S^{n-1}$
$$\int_{\R^n} |(\theta,\xi)|^{-1-q}\phi(\xi) d\xi = \frac{2\Gamma(-q)\cos(\pi q/2)}{\pi} \int_0^\infty t^{q}\hat{\phi}(t\theta) dt.$$
\end{lemma}

Our next lemma is a generalization of Theorem 3.18 from \cite{K1}, which was proved in \cite{GKS}.

\begin{lemma}
Let $K$ be an infinitely smooth origin-symmetric convex body in $\R^n,$ let $f$ be an even infinitely smooth function on $K,$ and let $q\in (-1,n-1).$ Then for every fixed $\xi\in S^{n-1}$
\begin{equation}\label{fourier}
(Rf(\xi,t))_t^{(q)}(0) 
\end{equation}
$$= \frac{\cos(\pi q/2)}{\pi} \left(|x|_2^{-n+q+1}\left(\int_0^{\frac{|x|_2}{\|x\|_K}} r^{n-q-2}f(r\frac{x}{|x|_2}) dr \right) \right)^\wedge_x(\xi).$$
\end{lemma}
\pf Let $-1<q<0.$ Then
$$(Rf(\xi,t))_t^{(q)}(0) = \frac 1{2\Gamma(-q)} 
\int_{-\infty}^{\infty} |t|^{-1-q}\left(\int_{K\cap \{x:\ (x,\xi)=t\}} f(x) dx\right) dt$$
$$=\frac 1{2\Gamma(-q)}  \int_{\R^n} |(x,\xi)|^{-1-q} f(x) \chi_K(x) dx $$
$$=\frac 1{2\Gamma(-q)}  \int_{S^{n-1}} |(\theta,\xi)|^{-1-q}\left( \int_0^{\|\theta\|_K^{-1}} r^{n-q-2} f(r\theta) dr\right) d\theta$$
Consider the latter as a homogeneous of degree $-1-q$ function of $\xi\in \R^n\setminus\{0\},$ apply it to an even test function $\phi$ and use Lemma \ref{lemma3.14}:
$$ \langle(Rf(\xi,t))_t^{(q)}(0), \phi \rangle$$
$$ = \dfrac{1}{2\Gamma(-q)} \int_{\R^n} \phi(\xi) \left(\int_{S^{n-1}} \left| ( \theta, \xi ) \right|^{-1 -q} \left( \int_0^{\|\theta\|_{K}^{-1}} r^{n-q-2} f(r\theta) dr\right) d\theta \right) d\xi $$
$$= \frac{\cos(q\pi/2)}{\pi} \int_{S^{n-1}} \left(\int_0^{\infty} t^q\hat{\phi}(t\theta) dt \right)\left(\int_0^{\|\theta\|_K^{-1}} r^{n-q-2} f(r\theta) dr \right) d\theta.$$
On the other hand, if we apply the function in the right-hand side of (\ref{fourier}) to the test function $\phi$ we get 
$$ \left\langle \dfrac{\cos(\pi q / 2)}{\pi} \left( |x|_2^{-n + q + 1} \left( \int_0^{\frac{|x|_2}{\|x\|_K}} r^{n-q-2} f(\frac{rx}{|x|_2}) dr\right) \right)_x^{\wedge} (\xi), \phi(\xi)\right\rangle $$
$$ =\dfrac{\cos(\pi q / 2)}{\pi} \int_{\R^n} |x|_2^{-n + q + 1} \left( \int_{0}^{\frac{|x|_2}{\|x\|_K}} r^{n-q-2} f(\frac{rx}{|x|_2}) dr\right) \hat{\phi}(x) dx$$
$$ =\dfrac{\cos(\pi q / 2)}{\pi} \int_{S^{n-1}} \left(\int_0^\infty t^q \hat{\phi} (t\theta) dt\right) \left( \int_0^{\|\theta\|_K^{-1}} r^{n-q-2} f(r\theta) dr\right) d\theta.$$
This proves the result for $-1<q<0.$ By a simple argument similar to that in Lemma 2.22 from \cite{K1}, one can see that both sides
of (\ref{fourier}) are analytic functions of $q$ in the domain $-1<\Re q <n-1.$
By analytic extension, (\ref{fourier}) holds for all $-1<q<n-1.$
\endpf

Let us compute the derivative of the Radon transform in the case where $f\equiv 1$ and $K=B_2^n,$ the unit Euclidean ball.
We denote by $\chi_K$ the indicator function of $K.$ Let $B_n=B_2^n/|B_2^n|^{\frac 1n}$ be the Euclidean ball of volume 1 in $\R^n.$
\begin{co} \label{eucl1} For $-1<q<n-1$ and every $\xi\in S^{n-1}$
$$(R(\chi_{B_2^n})(\xi,t))_t^{(q)}(0)=|B_2^n\cap \{x:(x,\xi)=t\}|_t^{(q)}(0)$$$$
=((1-t^2)^{\frac{n-1}2})^{(q)}(0) |B_2^{n-1}|= \frac{\cos(\frac{\pi q}2)}{\pi(n-q-1)} (|x|_2^{-n+q+1})^\wedge(\xi)$$
$$=\frac{2^{q+1} \pi^{\frac{n-2}2} \Gamma(\frac{q+1}2) \cos(\frac{\pi q}2)}{(n-q-1)\Gamma(\frac{n-q-1}2)}.$$
\end{co}

\pf Use (\ref{fourier}) with $f\equiv 1$ and $K=B_2^n$ and the formula for the Fourier transform of powers of the Euclidean norms from \cite{GS}: if $\lambda \in (-n,0),$ then
\begin{equation} \label{eucl}
(|x|_2^{\lambda})^\wedge (\xi) = \frac{2^{\lambda +n} \pi^{\frac{n}2} \Gamma(\frac{\lambda +n}2)}
{\Gamma(\frac{-\lambda}2)}|\xi|_2^{-\lambda-n}.
\end{equation}
\endpf
\bigbreak
Let $B_n=B_2^n/|B_2^n|^{\frac 1n}$ be the Euclidean ball of volume 1 in $\R^n,$ and $q\ge 0$ not an odd integer.
Then by Corollary \ref{eucl1} and since $|B_2^n|=\pi^{n/2}/\Gamma(\frac n2+1),$ $\Gamma(x+1)=x\Gamma(x)$,
and the $\Gamma$-function is log-convex (see \cite[Lemma 2.14]{K1}),
\begin{equation}\label{eucl-norm}
\frac 1{\cos(\frac{\pi q}2)} (R(\chi_{B_n})(\xi,t))_t^{(q)}(0)=\frac 1{\cos(\frac{\pi q}2)|B_2^n|^{\frac{n-q-1}n}}(R(\chi_{B_2^n})(\xi,t))_t^{(q)}(0)
\end{equation}
$$= \frac{2^{q} \pi^{\frac{n-2}2} \Gamma(\frac{q+1}2)\left(\Gamma(1+\frac n2)\right)^{\frac {n-q-1}n}}{\Gamma(\frac{n-q-1}2+1)\pi^{\frac{n-q-1}2}}\ge c^{q+1}(q+1)^{\frac{q+1}2},$$
where $c$ is an absolute constant.
\bigbreak

For $q>0,$ the radial $q$-sum of star bodies $K$ and $L$ in $\R^n$ is defined as a new star body $K\tilde{+}_qL$ whose radius 
in every direction $\xi\in S^{n-1}$ is given by
$$r^q_{K\tilde{+}_qL}(\xi)= r^q_{K}(\xi) + r^q_{L}(\xi),\qquad \forall \xi\in S^{n-1}.$$
The radial metric in the class of origin-symmetric star bodies is defined by
$$\rho(K,L)=\sup_{\xi\in S^{n-1}} |r_K(\xi)-r_L(\xi)|.$$

For the purposes of this article, we adopt the following definition of intersection bodies.
\begin{df}\label {bp-bodies}For $0<q<n,$ we define the class of {\it generalized $q$-intersection bodies} ${\cal{BP}}_q^n$ in $\R^n$ 
as the closure in the radial metric of radial $q$-sums of finite collections of origin-symmetric ellipsoids in $\R^n.$
\end{df}

When $q=1,$ by the result of Goodey and Weil \cite{GW}, we get the original class of intersection bodies ${\cal{I}}_n$ in $\R^n$ introduced by Lutwak \cite{Lu}. For integers $k,\ 1<k<n$, by the result of  Grinberg and Zhang \cite{GZ}, we get the classes 
of generalized $k$-intersection bodies introduced by Zhang \cite{Z2}. Note that the original definitions of Lutwak and Zhang are different; see \cite[Chapter 4]{K1} for details.

Another generalization of the concept of an intersection body was introduced 
in \cite{K-2000}. For $k=1,$ this is the original definition of Lutwak.

\begin{df}For an integer $k,\ 1\le k <n$ and star bodies $D,L$ in $\R^n,$
we say that $D$ is the $k$-intersection body of $L$ if 
$$|D\cap H^{\perp}|= |L\cap H|,\qquad \forall H\in Gr_{n-k}.$$
Taking the closure in the radial metric of the class of
$k$-intersection bodies of star bodies, we define the class of 
{$k$-intersection bodies}.
\end{df}

The latter definition is related to embeddings in $L_p$-spaces. By $L_p,\ p>0$ we mean the $L_p$-space of functions on $[0,1]$ with Lebesgue measure. It was proved in \cite{K-2000}  that an $n$-dimensional normed 
space embeds isometrically in $L_p,$ where $p>0$  
and $p$ is not an even integer, if and only if the Fourier transform in the sense of Schwartz distributions
of the function $\Gamma(-p/2)\|\cdot\|^p$ is a non-negative distribution outside of the origin in $\R^n.$
The concept of embedding of finite dimensional normed spaces in $L_p$ with negative $p$
was introduced in \cite{K-CMB, K-2000}, as an 
extension of embedding into $L_p$ with $p>0,$ as follows.

\begin{df} For a star body $D$ in $\R^n,$ we say that $(\R^n,\|\cdot\|_D)$ embeds in $L_{-p},\ 0<p<n$ 
if the function $\|\cdot\|_D^{-p}$ represents a positive definite distribution.
We denote the class of such star bodies by $L_{-p}^n.$
\end{df}

A connection between $k$-intersection bodies and embedding in $L_p$ with
negative $p$ was found in \cite{K-2000}.
\begin{theorem} \label{connection} Let $1\le k < n.$ The following are equivalent:
\item{(i)} An origin symmetric star body $D$ in $\R^n$ is a $k$-intersection
body;
\item{(ii)} The space $(\R^n,\|\cdot\|_D)$ embeds in $L_{-k}.$
\end{theorem} 

It was proved in \cite{K-2000} (see also \cite{M2, KY}) that for any integer $k,\ 1\le k <n$ every generalized $k$-intersection body is a $k$-intersection body. We need an extension of this fact to non-integers, as follows.

\begin{lemma} \label{incl} For every $0<q<n,$ we have ${\cal{BP}}_q^n\subset L_{-q}^n.$
\end{lemma} 
\pf The powers of the Euclidean norm $|x|_2^{-q},\ 0<q<n$ represent positive definite distributions in $\R^n;$ see formula (\ref{eucl}). 
Because of the connection between the Fourier transform of distributions and linear transformations, the norms generated by origin-symmetric ellipsoids in $\R^n$ have the same property, and, therefore, radial $q$-sums of ellipsoids in $\R^n$ belong to
the class $L_{-q}^n.$ The fact that positive definiteness is preserved under limits in the radial metric follows from
\cite[Lemma 3.11]{K1}. \endpf

\section{Proofs of main results}
 
 \noindent{\bf Proof of Theorem \ref{bp}.} Let $D$ be such that $|D|^{1/n}\le (1+\delta) d_{ovr}(K,L^n_{-1-q}),\ K\subset D,\ D\in I_q,\ \delta>0.$ By approximation, we can assume that $D$ is infinitely smooth; see \cite[Lemma 4.10]{K1}.
 Then $(\|x\|_D^{-1-q})^\wedge$ is a positive function on the sphere. By (\ref{fourier}) for every $\xi\in S^{n-1}$
 $$\left(|x|_2^{-n+q+1}\left(\int_0^{\frac{|x|_2}{\|x\|_K}} r^{n-q-2}f(r\frac{x}{|x|_2}) dx \right) \right)^\wedge_x(\xi) $$$$\le
\left(|x|_2^{-n+q+1}\left(\int_0^{\frac{|x|_2}{\|x\|_L}} r^{n-q-2}g(r\frac{x}{|x|_2}) dx \right) \right)^\wedge_x(\xi).$$
Multiplying both sides by $(\|x\|_D^{-1-q})^\wedge(\xi)$, integrating over the sphere and using Parseval's formula on the sphere
(\ref{parseval}) we get
$$\int_{S^{n-1}} \|\theta\|_D^{-1-q}\left( \int_0^{\|\theta\|_K^{-1}} r^{n-q-2} f(r\theta) dr\right)$$
$$\le \int_{S^{n-1}} \|\theta\|_D^{-1-q}\left( \int_0^{\|\theta\|_L^{-1}} r^{n-q-2} g(r\theta) dr\right),$$
or
\begin{equation}\label{sect12}
\int_K \|x\|_D^{-1-q} f(x) dx \le \int_L \|x\|_D^{-1-q} g(x) dx.
\end{equation}
Since $K\subset D,$ we have $1\ge \|x\|_K\ge \|x\|_D$ for every $x\in K.$ Therefore,  
$$\int_K \|x\|_D^{-1-q}f(x) dx \ge \int_K \|x\|_K^{-1-q}f(x) dx \ge \int_K f(x) dx.$$

On the other hand, by the Lemma from section 2.1 from Milman-Pajor
\cite[p.76]{MP}, 
$$
\left(\frac{\int_{L}\|x\|_D^{-1-q} {g}(x) dx}{ \int_{D}\|x\|_D^{-1-q} dx} \right)^{1/(n-q-1)}  \le \left(\frac{\int_{L}{g}(x) dx \large}{  \int_{D} dx} \right)^{1/n}.
$$
Since $\int_D\|x\|_D^{-1-q} dx =\frac{n}{n-q-1} |D|,$  (see \cite{MP}) we can estimate the right-hand side of
(\ref{sect12}) by
$$\int_L \|x\|_D^{-1-q} g(x) dx \le \frac n{n-q-1}  \left(\int_{L}{g}(x) dx\right)^{\frac {n-q-1}n} |D|^{\frac {q+1}n}.$$
Combining these estimates with the definition of $D$ and sending $\delta$ to zero, we get the result.
\endpf

 \noindent{\bf Proof of Theorem \ref{slicing}.} Consider a number $\e>0$ such that for every $\xi\in S^{n-1}$
 $$\frac 1{\cos(\pi q/2)}(R(f\vert_K)(\xi,t))_t^{(q)}(0)\le \frac {\e}{\cos(\pi q/2)}R(\chi_{B_2^n})(\xi,t))_t^{(q)}(0).$$
By (\ref{fourier}) for every $\xi\in S^{n-1}$
 $$\left(|x|_2^{-n+q+1}\left(\int_0^{\frac{|x|_2}{\|x\|_K}} r^{n-q-2}f(r\frac{x}{|x|_2}) dx \right) \right)^\wedge_x(\xi) \le
\frac{\e}{n-q-1}\left(|x|_2^{-n+q+1}\right)^\wedge(\xi)$$
Let $D$ be such that $|D|^{1/n}\le (1+\delta) d_{\rm ovr}(K,L^n_{-1-q}),\ K\subset D,\ D\in I_q,\ \delta>0.$ By approximation, we can assume that $D$ is infinitely smooth.
 Then $(\|x\|_D^{-1-q})^\wedge$ is a positive function on the sphere.. 
Multiplying both sides of the latter inequality by $(\|x\|_D^{-1-q})^\wedge(\xi)$, integrating over the sphere and using Parseval's formula on the sphere
we get
\begin{equation}\label{sect13}
\int_K \|x\|_D^{-1-q} f(x) dx \le \frac{\e}{n-q-1}\int_{S^{n-1}} \|x\|_D^{-1-q} dx.
\end{equation}
Since $K\subset D,$ we have 
\begin{equation}\label{eq14}
\int_K \|x\|_D^{-1-q} f(x) dx \ge \int_K \|x\|_K^{-1-q} f(x) dx \ge  \int_K f.
\end{equation}
On the other hand, by H\"older's inequality and the polar formula for volume,
\begin{equation}\label{eq15}
\int_{S^{n-1}} \|x\|_D^{-1-q} dx \le 
|S^{n-1}|^{\frac{n-q-1}n} n^{\frac{q+1}n} |D|^{\frac{q+1}n}.
\end{equation}
Put 
$$\e=\max_{\xi\in S^{n-1}} \frac{\frac 1{\cos(\pi q/2)}(R(f\vert_K)(\xi,t))_t^{(q)}(0)}{\frac {1}{\cos(\pi q/2)}R(\chi_{B_2^n})(\xi,t))_t^{(q)}(0)},$$
then (\ref{sect13}), (\ref{eq14}), (\ref{eq15}) and the choice of $D$ (send $\delta$ to zero) imply
$$\int_K f\le c(n,q) |D|^{\frac{q+1}n} \max_{\xi\in S^{n-1}} \frac 1{\cos(\pi q/2)}(R(f\vert_K)(\xi,t))_t^{(q)}(0)$$
$$\le c(n,q) \left(d_{\rm ovr}(K,L^n_{-1-q})\right)^{q+1} |K|^{\frac{q+1}n} \max_{\xi\in S^{n-1}} \frac 1{\cos(\pi q/2)}(R(f\vert_K)(\xi,t))_t^{(q)}(0),$$
 where by Corollary \ref{eucl1}, the property $\Gamma(x+1)=x\Gamma(x)$ of the $\Gamma$-function and the expression for the surface area of the unit Euclidean sphere (see for example \cite[Corollary 2.20]{K1}),
$$c(n,q)= \frac{\pi \Gamma(\frac{n-q-1}2) n^{\frac{q+1}n} 2^{\frac{n-q-1}n} \pi^{\frac{n-q-1}2}}
{2^{q+1}\pi^{\frac n2} \Gamma(\frac{q+1}2) \left(\Gamma(\frac n2)\right)^{\frac{n-q-1}n}}$$
$$= \frac{n\ \Gamma(\frac{n-q-1}2+1)}{2^q \pi^{\frac{q-1}2} \Gamma(\frac{q+1}2) (n-q-1) 
\left(\Gamma(\frac n2 +1)\right)^{\frac{n-q-1}n}}\le  \frac{n}{2^q \pi^{\frac{q-1}2} \Gamma(\frac{q+1}2) (n-q-1)},$$
where in the last step we used the log-convexity of the $\Gamma$-function (see \cite[Lemma 2.14]{K1}) to estimate
$$\frac{\Gamma(\frac{n-q-1}2+1)}{\left(\Gamma(\frac n2 +1)\right)^{\frac{n-q-1}n}}\le 1.$$
\endpf
 
To prove Theorem \ref{fract-der} we need the estimate for the outer volume ratio distance from \cite[Theorem 1.1]{KPZ}. Note that this estimate was proved in \cite{KPZ} for integers $q,$ but the proof remains exactly the same for non-integers. Also note that a mistake in the proof in \cite{KPZ} was corrected in \cite[Section 5]{K2}.

\begin{pr} \label{kpz} (\cite{KPZ}) For every $q\in [1,n-1]$ and every origin-symmetric convex body $K$ in $\R^n$
$$d_{\rm ovr}(K, {\cal{BP}}_q^n)\le C\sqrt{\frac{n\log^3(\frac{ne}{q})}{q}},$$
where $C$ is an absolute constant, and 
$$d_{{\rm {ovr}}}(K, {\cal{BP}}_q^n) = \inf \left\{ \left( \frac {|D|}{|K|}\right)^{1/n}:\ K\subset D,\ D\in {\cal{BP}}_q^n \right\}.$$ 
\end{pr}

\noindent{\bf Proof of Theorem \ref{fract-der}.} By Lemma \ref{incl} and Proposition \ref{kpz},
\begin{equation}\label{general}
d_{\rm ovr}(K,L^n_{-1-q})\le d_{\rm ovr}(K, {\cal{BP}}_{q+1}^n)\le C\sqrt{\frac{n\log^3(\frac{ne}{q+1})}{q+1}}.
\end{equation}
Now the result follows from Theorem \ref{slicing}. Note that 
 $$\frac n{n-q-1} \le e^{\frac{q+1}{n-q-1}} \le e^{q+1}.$$
 \endpf
 
 We conclude this note with several remarks about the classes $L_{-p}^n$ and the outer volume ratio distance 
 $d_{\rm ovr}(K,L^n_{-p}).$ First, every class $L_{-p},\ 0<p<n$ contains the unit balls of all $n$-dimensional subspaces of $L_r$ with $0<r\le 2;$ see \cite{K5} and \cite[Theorem 6.17]{K1}. In particular all of these classes contain the $\ell_1^n$-ball
 $$B_1^n=\{x\in \R^n: |x_1|+...+|x_n|\le 1\}.$$
 This implies that the distance  $d_{\rm ovr}(K,L^n_{-p})\le e$ for every unconditional convex body $K$ and every $0<p<n;$
 see \cite{K2}. Since all the classes $L_{-p}$ contain origin-symmetric ellipsoids, by a result from \cite{M1} (see also \cite{KP}),
 the distance $d_{\rm ovr}(K,L^n_{-p})\le c\sqrt{r},$ where $K$ is the unit ball of an $n$-dimensional subspace of $L_r$ with $r>2,$ and $c$ is an absolute constant. It was proved in \cite{M2} that the class $L_{-m}^n$ contains $L_{-k}^n$ when $m,k$ are integers and $m$ is divisible by $k.$ This inclusion for arbitrary $0<k<m<n$ has not been established in general. It is only known that if $p>s,$ then $L_{-p}^n$ contains a body that is not in $L_{-s}^n;$ see \cite{KaZ, Y}.  When If $n-3\le p<n,$ the class
 $L_{-p}^n$ contains all origin-symmetric convex bodies in $\R^n;$ see \cite[Corollary 4.9]{K1}. More results about embeddings in $L_{-p}$ can be found in \cite{KaK, KPZ, K6, M1, M2}, \cite[Chapter 6]{K1} and \cite{KY}. It would be interesting to know whether the estimate (\ref{general}) is sharp.

\end{document}